\documentclass{amsart}
\usepackage{latexsym}   
\usepackage{amssymb}    
\usepackage{amsmath}    
\usepackage{amsthm}     
\usepackage[all]{xy}
\usepackage{hyperref}

\newcommand{\pa}{\partial}

\newtheorem{theorem}{Theorem}[section]
\newtheorem{lemma}[theorem]{Lemma}

\theoremstyle{definition}

\theoremstyle{remark}
\newtheorem{remark}[theorem]{Remark}

\numberwithin{equation}{section}

\begin{document}

\title{Archimedes' balance and Bianchi's B\"{a}cklund transformation for quadrics}

\author{Ion I. Dinc\u{a}}
\address{Faculty of Mathematics and Informatics,
University of Bucharest,  14 Academiei Str., 010014, Bucharest,
Romania}
 \email{dinca@gta.math.unibuc.ro}
\thanks{Supported by the University of Notre Dame du Lac}



\begin{abstract}
We establish a link between Archimedes' method of integration for
calculating areas, volumes and centers of mass of segments of
parabolas and quadrics of revolution by factorization via the
moments of a balance and an integration technique for a particular
integrable system, namely Bianchi's B\"{a}cklund transformation
for quadrics.
\end{abstract}

\maketitle

\section*{Introduction}
This paper is organized by first stating the relevant results of
Archimedes and Bianchi as they originally appeared, then having a
discussion on the notions of Bianchi's result and explaining them
in terms of current definitions, providing short motivation and
proof for real ruled quadrics for Bianchi's result and then
explaining the link between the results of the two authors.

\section{Archimedes' and Bianchi's results}

In {\it The Method} as it appears in \cite{M} Archimedes claims:

\begin{theorem}[The Method]
'... certain things first became clear to me by a mechanical
method, although they had to be proved by geometry afterwards
because their investigation by the said method did not furnish an
actual proof. But it is of course easier, when we have previously
acquired, by the method, some knowledge of the questions, to
supply the proof than it is to find it without any previous
knowledge.'
\end{theorem}

The main theorem of Bianchi's theory of deformations of surfaces
applicable to quadrics (which proves the existence of the
B\"{a}cklund transformation, its inversion and of the
applicability correspondence provided by the Ivory affinity)
roughly states:
\begin{theorem}[Theorem I]

Every surface $x^0\subset \mathbb{C}^3$ applicable to a surface
$x_0^0\subseteq x_0$ ($x_0$ being a quadric) appears as a focal
surface of a $2$-dimensional family of Weingarten congruences,
whose other focal surfaces $x^1=B_{z}(x^0)$ (called B\"{a}cklund
transforms of $x^0$) are applicable, via the Ivory affinity, to
surfaces $x_0^1$ in the same quadric $x_0$. The determination of
these surfaces requires the integration of a family of Riccatti
equations depending on the parameter $z$ (ignore for simplicity
the dependence on the initial value of the Ricatti equation in the
notation $B_z$). Moreover, if we compose the inverse of the rigid
motion provided by the Ivory affinity with the rolling of $x_0^0$
on $x^0$, then we obtain the rolling of $x_0^1$ on $x^1$ and $x^0$
reveals itself as a $B_z$ transform of $x^1$.
\end{theorem}

\section{Discussion on the notions appearing in Bianchi's result}

\begin{remark}
The use of imaginaries (when one complexifies both the surface and
the surrounding Euclidean space) is important because it is Lie's
interpretation of the B\"{a}cklund transformation for constant
Gau\ss\ curvature $-1$ surfaces on (imaginary) confocal
pseudo-spheres\footnote{ From Bianchi's quotation it is unclear to
us if this point of view is due entirely to Lie or if Lie observed
the collapsing of leaves of $2$-dimensional integrable
distributions of facets to curves and points and further Bianchi
made the remark for the pseudo-sphere.} that was the tool that
allowed Bianchi to prove his result; for this reason we chose to
state it in a complex setting.
\end{remark}

Except for Lie's influence we shall only work with objects
immersed in the Euclidean space

$$(\mathbb{R}^3,<\cdot,\cdot>),\ <x,y>:=x^Ty,\ |x|^2:=x^Tx\
\mathrm{for}\ x,y\in\mathbb{R}^3.$$

The standard basis $\{e_1,e_2,e_3\}$ satisfies
$e_i^Te_j=\delta_{ij}$.

In this setting {\it 'applicable'} surfaces means just 'isometric
surfaces' and 'applicability correspondence' means just 'isometric
correspondence' (local diffeomorphism).

A {\it 'Weingarten congruence'} is a $2$-dimensional family of
lines on whose two focal surfaces the asymptotic coordinates
correspond (equivalently the second fundamental forms of the two
focal surfaces are proportional).

\begin{remark}
Note that although the correspondence provided by the Weingarten
congruence is not the isometric one, a Weingarten congruence is
the tool best suited to attack the isometric deformation problem
by means of transformation, since it provides correspondence of
the characteristics of the isometric deformation problem
(according to Darboux these are the asymptotic coordinates) and it
is directly linked to the infinitesimal isometric deformation
problem (Darboux proved that infinitesimal isometries generate
Weingarten congruences and Guichard proved the converse).
\end{remark}

Since we want both the {\it seed} $x^0$ and the {\it leaves}
$x^1=B_z(x^0)$ to be real surfaces isometric to pieces of real
quadrics, the quadric $x_0$ and its {\it confocal} (with same
foci) one $x_z$ must be real doubly ruled, so $x_0, x_z$ are
either hyperboloids with one sheet or hyperbolic paraboloids.

The Ivory affinity between confocal quadrics is a natural affine
correspondence between confocal quadrics and having good metric
properties; thus we have the hyperboloid with one sheet

\begin{eqnarray}\label{eq:Qhyp}
x_z(u,v):=\sqrt{a_1-z}\frac{1-uv}{u-v}e_1+\sqrt{z-a_2}\frac{1+uv}{u-v}e_2+
\sqrt{a_3-z}\frac{u+v}{u-v}e_3,\nonumber\\ a_2<0,z<a_1,a_3;\
u,v\in \mathbb{R}\cup\{\infty\},\ u\neq v,
\end{eqnarray}
when the Ivory affinity is given by

$$x_z(u,v)=\sqrt{I_3-zA}\
x_0(u,v),\ A:=\mathrm{diag}[a_1^{-1}\ \ a_2^{-1}\ \ a_3^{-1}]$$

and the hyperbolic paraboloid
\begin{eqnarray}\label{eq:Qpar}
x_z(u,v):=\sqrt{a_1-z}(u+v)e_1+\sqrt{z-a_2}(u-v)e_2+(2uv+\frac{z}{2})e_3,\nonumber\\
a_2<0,z<a_1,\ u,v\in\mathbb{R},
\end{eqnarray}
when the Ivory affinity is given by $$x_z(u,v)=\sqrt{I_3-zA}\
x_0(u,v)+\frac{z}{2}e_3,\ A:=\mathrm{diag}[a_1^{-1}\ \ a_2^{-1}\ \
0].$$

For $D\subseteq \mathbb{R}^2$ domain two isometric surfaces $x_0,\
x:D\rightarrow \mathbb{R}^3,\ |dx_0|^2=|dx|^2$ can be rolled one
onto the other: $(x,dx)=(R,t)(x_0,dx_0):=(Rx_0+t,Rdx_0),\\
(R,t):D\rightarrow \mathbf{O}_3(\mathbb{R})\ltimes\mathbb{R}^3$
being a surface in the space of rigid motions (it degenerates to a
curve if $x_0,\ x$ are ruled with isometric correspondence of
rulings, when the rolling takes place in a $1$-dimensional
fashion, or to a point if $x_0,\ x$ are rigidly isometric) such
that at any instant they meet tangentially and with same
differential at the tangency point:

\begin{eqnarray}\label{eq:dx}
dx=Rdx_0.
\end{eqnarray}
Conversely, if $x_0$ can be rolled on $x$ (that is we have
(\ref{eq:dx})), then $x_0$ and $x$ are isometric.

For $(u,v)$ parametrization on $D$ and $N_0,\ N$ Gau\ss\ maps
respectively of $x_0,\ x$ we have $R[\pa_ux_0\ \ \pa_vx_0\ \
N_0]=[\pa_ux\ \ \pa_vx\ \ \epsilon N],\ \epsilon:=\pm 1$, so the
rotation $R$  of the rolling is uniquely defined (modulo the
indeterminacy $\epsilon$) by $$R:=[\pa_ux\ \ \pa_vx\ \ \epsilon
N][\pa_ux_0\ \ \pa_vx_0\ \ N_0]^{-1};$$ the translation $t$ is
then given by $t:=x-Rx_0$. The indeterminacy $\epsilon$ decides
wether $R$ is special orthogonal or not (we have
$N=\frac{\pa_ux\times\pa_vx}{|\pa_ux\times
\pa_vx|}=\frac{R\pa_ux_0\times R\pa_vx_0}{|R\pa_ux_0\times
R\pa_vx_0|}=\det(R)(R^T)^{-1}\frac{\pa_ux_0\times
\pa_vx_0}{|\pa_ux_0\times \pa_vx_0|}=\det (R)RN_0$) and it has a
simple geometric explanation: $x_0$ can be rolled on either side
of $x$. It is an immediate consequence of (\ref{eq:dx}), since
(\ref{eq:dx}) involves only information about the tangent bundle,
so symmetries of the normal bundles (reflections in surfaces) are
allowed.

\begin{remark}\label{rm:refls}
Although Bianchi was aware of this indeterminacy, its importance
seems to have escaped his attention; we shall see later that this
indeterminacy provides a simple geometric explanation of an
indeterminacy appearing in the rigid motion provided by the Ivory
affinity (choice of ruling), which in turn encodes all necessary
algebraic information needed to prove Bianchi's result.
\end{remark}

\begin{remark}
We shall use the notation $d$ for exterior (antisymmetric)
derivative; thus $d^2=0$.
\end{remark}

For $\omega_1,\ \omega_2\ \mathbb{R}^3$-valued $1$-forms on $D$
and $a,b\in\mathbb{\mathbb{R}}^3$, we have $a^T\omega_1\wedge
b^T\omega_2=((a\times
b)\times\omega_1+b^T\omega_1a)^T\wedge\omega_2=(a\times
b)^T\omega_1\times\wedge\omega_2+b^T\omega_1\wedge a^T\omega_2$;
in particular
\begin{eqnarray}\label{eq:oma}
a^T\omega\wedge b^T\omega=\frac{1}{2}(a\times
b)^T\omega\times\wedge\omega.
\end{eqnarray}

Since both $\times$ and $\wedge$ are skew-symmetric, we have
$2\omega_1\times\wedge\omega_2=
\omega_1\times\omega_2+\omega_2\times\omega_1=2\omega_2\times\wedge\omega_1$.

Applying $d$ to (\ref{eq:dx}) and then multiplying with $R^{-1}$
we get $R^{-1}dR\wedge dx_0=0$; since $R^{-1}dR$ is
skew-symmetric, we have $dx_0^TR^{-1}dRdx_0=0$.

For $a\in\mathbb{R}^3$ we have
$R^{-1}dRa=R^{-1}dR(a^{\bot}+a^{\top})=a^TN_0R^{-1}dRN_0-\\a^TR^{-1}dRN_0N_0=\omega_0\times
a,\ \omega_0:=N_0\times R^{-1}dRN_0$, so under the identification
$(\mathbf{o}_3(\mathbb{R}),[,])\simeq(\mathbb{R}^3,\times)$ we
have $R^{-1}dR\simeq\omega_0$. Imposing the compatibility
condition $d$ to $R^{-1}dR$ and to (\ref{eq:dx}) we get
\begin{eqnarray}\label{eq:Qom}
d\omega_0+\frac{1}{2}\omega_0\times\wedge\omega_0=0,\
\omega_0\times\wedge dx_0=0
\end{eqnarray}

and thus $\omega_0$ is a flat connection form in $Tx_0$ (it
encodes the difference between the Gau\ss\ -Codazzi-Mainardi
equations for $x, x_0$).

Finally the rigid motion provided by the Ivory affinity exists due
to two results of Ivory's and Bianchi's on confocal quadrics (some
of the properties used by Bianchi may have already been folklore
by that time; for example the change in angles between rulings on
confocal hyperbolic paraboloids while preserving their lengths was
known to Henrici when he constructed the articulated hyperbolic
paraboloid).

\begin{theorem}[Ivory]
The orthogonal trajectory of a point on $x_z$ (as $z$ varies) is a
conic and the correspondence $x_0\rightarrow x_z$ thus established
is affine.

This affine transformation (henceforth called the Ivory affinity)
preserves the lengths of segments between confocal quadrics: with
$V_0^1:=x_z^1-x_0^0,\ V_1^0:=x_z^0-x_0^1$ we have
$|V_0^1|^2=|V_1^0|^2$ for pairs of points $(x_0^0,x_z^0),\
(x_0^1,x_z^1)$ corresponding on $(x_0,x_z)$ under the Ivory
affinity.
\end{theorem}

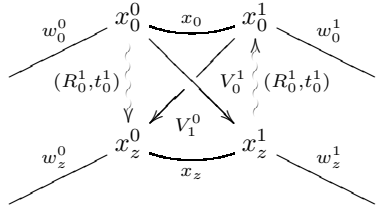
\begin{center}$\xymatrix@!0{&&x_0^0\ar@{-}[drdr]\ar@/_/@{-}[rr]^{x_0}\ar@{~>}[dd]_{(R_0^1,t_0^1)}&&
x_0^1\ar@{<~}[dd]^{(R_0^1,t_0^1)}&&\\
\ar@{-}[urr]^{w_0^0}&&&\ar[dl]^>>>>{V_1^0}\ar[dr]^>>>>>{V_0^1}&&&\ar@{-}[ull]_{w_0^1}\\&&
x_z^0\ar@{-}'[ur][urur]\ar@/_/@{-}[rr]_{x_z}&&x_z^1&\\\ar@{-}[urr]^{w_z^0}&&&
&&&\ar@{-}[ull]_{w_z^1}}$
\end{center}
\begin{theorem}[Bianchi I]\label{thm:B1}
If we have the rulings $w_0^0,\ w_0^1$ at the points $x_0^0,\
x_0^1\in x_0$ and by use of the Ivory affinity we get the rulings
$w_z^0,\ w_z^1$ at the points $x_z^0,\ x_z^1\in x_z$, then
$[V_0^1\ \ w_0^0\ \ w_z^1]^T[V_0^1\ \ w_0^0\ \ w_z^1]=[-V_1^0\ \
w_z^0\ \ w_0^1]^T[-V_1^0\ \ w_z^0\ \ w_0^1]$, so there exists a
rigid motion
$(R_0^1,t_0^1)\in\mathbf{O}_3(\mathbb{R})\ltimes\mathbb{R}^3$ with
\begin{eqnarray}\label{eq:Qivo}
(R_0^1,t_0^1)(x_0^0,x_z^1,w_0^0,w_z^1)=(x_z^0,x_0^1,w_z^0,w_0^1).
\end{eqnarray}
Moreover
$(V_0^1)^T\partial_z|_{z=0}x_z^0=(V_1^0)^T\partial_z|_{z=0}x_z^1$,
so the Ivory affinity has a nice projective property: the symmetry
of the tangency configuration
\begin{eqnarray}\label{eq:Qtc}
x_z^1\in T_{x_0^0}x_0\Leftrightarrow x_z^0\in T_{x_0^1}x_0.
\end{eqnarray}
\end{theorem}
\begin{remark}
The action of the rigid motion provided by the Ivory affinity
resembles a balance; in fact Bianchi uses moments and angles of
pairs of lines to explain it.
\end{remark}

\section{Motivation and proof of Bianchi's main theorem}

\subsection{The Bianchi-Lie ansatz}

The transformation originally constructed by

B\"{a}cklund in 1883 states that if a constant Gau\ss\ curvature
$-1$ seed $x$ in general position and an angle
$0<\theta<\frac{\pi}{2}$ are given, then the $3$-dimensional
distribution formed by {\it facets} (pairs of points and planes
passing through those points; it is the infinitesimal version of a
surface) with centers on circles of radius $\sin\theta$ in tangent
planes of $x$ (the circles are themselves centered at the origins
of tangent planes), of inclination $\theta$ to these and passing
through the origin of these is integrable; moreover the leaves are
constant Gau\ss\ curvature $-1$ surfaces and their determination
requires the integration of a Ricatti equation.

The B\"{a}cklund transformation when $\theta=\frac{\pi}{2}$ was
constructed even earlier (1879, upon some results of Ribaucour
from 1870) by Bianchi in his PhD thesis and named the
complementary transformation, but B\"{a}cklund's merit is the
introduction of the spectral parameter.

Lie considered the natural question: because the B\"{a}cklund
transformation is of a general nature (independent of the shape of
$x$), it must exist (at least as a limiting case) also in the case
when $x$ in not in general position, but it actually coincides
with the pseudo-sphere. In this case the $1$-dimensional family of
non-degenerate leaves (surfaces) degenerates to a $1$-dimensional
family of degenerate leaves (isotropic rulings on confocal
pseudo-sphere) and thus the true nature of the B\"{a}cklund
transformation is revealed at the static level of confocal
pseudo-spheres.

Thus we consider together with Bianchi (\cite{B2},\S\ 374)

\begin{theorem}[Lie's inverting point of view]
The tangent planes to the unit pseudo-sphere $x_0$ cut a confocal
pseudo-sphere $x_z$ along circles, thus highlighting a circle in
each tangent plane of $x_0$. Each point of the circle, the segment
joining it with the origin of the tangent plane and one of the
(imaginary) rulings on $x_z$ passing through that point determine
a facet. We have thus highlighted a $3$-dimensional integrable
distribution of facets: its leaves are the ruling families on
$x_z$. If we roll the distribution while rolling $x_0$ on an
isometric surface $x$ (called {\it seed}), it turns out that the
integrability condition of the rolled distribution is always
satisfied (we have complete integrability), so the integrability
of the rolled distribution does not depend on the shape of the
seed. The rolled distribution is obtained as follows: each facet
of the original distribution corresponds to a point on $x_0$; we
act on that facet with the rigid motion of the rolling
corresponding to the highlighted point of $x_0$ in order to obtain
the corresponding facet of the rolled distribution. The leaves of
the rolled distribution (called the B\"{a}cklund transforms of
$x$, denoted $B_z(x)$ and whose determination requires the
integration of a Ricatti equation) are isometric to the
pseudo-sphere. Moreover the seed and any leaf are the focal
surfaces of a Weingarten congruence, so the inversion of the
B\"{a}cklund transformation has a simple geometric explanation
(the seed and leaf exchange places).
\end{theorem}

Lie's inverting point of view allows us to call the B\"{a}cklund
transformation of constant Gau\ss\ curvature $-1$ surfaces {\it
B\"{a}cklund transformation of the pseudo-sphere}.

\begin{remark}\label{rm:refl}
Note that the B\"{a}cklund transformation of the pseudo-sphere
comes in two flavors, as the facets may reflect in the tangent
plane upon which their centers lie, but the complementary
transformation comes only in one flavor; this corresponds to a
non-degenerate confocal pseudo-sphere having two distinct families
of imaginary rulings and respectively to the two families of
rulings degenerating to a single one on the light cone, a singular
confocal pseudo-sphere.
\end{remark}

\begin{remark}
Lie's ansatz is the one susceptible for generalization, since it
provides a geometric explanation of the dependence of the
B\"{a}cklund transformation on the spectral parameter $z$. Thus if
in Lie's interpretation one replaces {\it 'pseudo-sphere'} with
{\it 'quadric'} and {\it 'circle'} with {\it 'conic'}, then one
gets Bianchi's result except for the Ivory affinity influence.
\end{remark}

While looking for the isometric correspondence Bianchi rolled back
the seed on the original quadric, in which case the facets of the
rolled distribution return to their original location on the
confocal quadric $x_z$. Thus if one assumes the isometric
correspondence to be valid and of a general nature, it must be
independent of the shape of the seed and the answer should be
found on the confocal family: the Ivory affinity provides a
natural correspondence between confocal quadrics and proving that
it provides the isometric correspondence remained a matter of
computations.

We are actually able at this point to prove the isometric
correspondence provided by the Ivory affinity and the inversion of
the B\"{a}cklund transformation for quadrics, assuming {\bf
Theorem \ref{thm:B1}} and that the complete integrability of the
rolled distribution is checked (we need the leaf $x^1$ to exist):
if we roll the seed $x^0$ on $x_0^0$, then the tangent plane of
the leaf $x^1$ corresponding to the point of tangency of the
rolled $x_0^0$ and the seed $x^0$ will be applied to the facet
centered at $x_z^1$ and spanned by $V_0^1,\ \pa_{u_1}x_z^1$;
further applying the rigid motion $(R_0^1,t_0^1)$ provided by the
Ivory affinity it will be applied to $T_{x_0^1}x_0$. In this
process $\pa_{u_1}x^1$ is taken to $\pa_{u_1}x_z^1$ and further to
$\pa_{u_1}x_0^1$, so actually $(x^1,dx^1)$ is taken to
$(x_0^1,dx_0^1)$; moreover because of the symmetry of the tangency
configuration the seed becomes leaf and the leaf becomes seed.

\begin{remark}
This geometric argument was the one used by Bianchi to prove the
inversion of the B\"{a}cklund transformation, but it seems that he
preferred the security of an analytic confirmation to the power of
his geometric arguments for the isometric correspondence provided
by the Ivory affinity.
\end{remark}
\subsection{Proof of Bianchi I} (which includes also the Ivory theorem).

Note that with $B:=0,\ C:=-1$ for (\ref{eq:Qhyp}), $B:=-e_3,\
C:=0$ for (\ref{eq:Qpar}) and $R_z:=I_3-zA$ both confocal families
(\ref{eq:Qhyp}) and (\ref{eq:Qpar}) (in fact all confocal families
of quadrics) can be implicitly defined by
$$\begin{bmatrix}x_z\\1\end{bmatrix}^T(\begin{bmatrix}A&B\\B^T&C\end{bmatrix}^{-1}
-z\begin{bmatrix}I_3&0\\0&0\end{bmatrix})^{-1}\begin{bmatrix}x_z\\1\end{bmatrix}=0,$$
equivalently

$$\begin{bmatrix}x_z\\1\end{bmatrix}^T\begin{bmatrix}AR_z^{-1}&R_z^{-1}B\\
B^TR_z^{-1}&C+zB^TR_z^{-1}B\end{bmatrix}\begin{bmatrix}x_z\\1\end{bmatrix}=0.$$

From the first definition one can see the metric-projective
definition of the family of confocal quadrics: a pencil behavior
and Cayley's absolute
$\begin{bmatrix}x\\0\end{bmatrix}^T\begin{bmatrix}I_3&0\\0&0\end{bmatrix}
\begin{bmatrix}x\\0\end{bmatrix}=0,\ x\neq 0$ in the plane at $\infty$ (which
encodes the Euclidean structure on $\mathbb{R}^3$).

With   $C(z):=(-\frac{1}{2}\int_0^z(\sqrt{R_w})^{-1}dw)B(=0$ for
(\ref{eq:Qhyp}) and $=\frac{z}{2}e_3$ for (\ref{eq:Qpar})) we also
have an unifying formula for the Ivory affinity, namely
$x_z=\sqrt{R_z}x_0+C(z)$ and finally it is convenient to work with
the normal field $\hat N_z:=-2\partial_zx_z$ instead of with the
unit normal $N_z$ (note
$AC(z)+(I_3-\sqrt{R_z})B=0=(I_3+\sqrt{R_z})C(z)+zB$, since both
are $0$ for $z=0$ and do not depend on $z$).

Ivory  becomes:
\begin{itemize}
\item
$|V_0^1|^2=|x_0^0+x_0^1-C(z)|^2-2(x_0^0)^T(I_n+\sqrt{R_z})x_0^1+zC=|V_1^0|^2$;
\end{itemize}
Bianchi I becomes:  if $w_0^TAw_0=w_0^T\hat N_0=0,\
w_z=\sqrt{R_z}w_0$, etc, then:
\begin{itemize}
\item for lengths of rulings:
$w_z^Tw_z=|w_0|^2-zw_0^TAw_0=|w_0|^2$; \item for angles between
segments and rulings: $(V_0^1)^Tw_0^0+(V_1^0)^Tw_z^0=\\-z(\hat
N_0^0)^Tw_0^0=0$; \item for angles between rulings:
$(w_0^0)^Tw_z^1=(w_0^0)^T\sqrt{R_z}w_0^1=(w_z^0)^Tw_0^1$; \item
for the symmetry of the tangency configuration: $(V_0^1)^T\hat
N_0^0
=(x_0^0)^TA\sqrt{R_z}x_0^1-B^T(x_z^0+x_z^1-C(z))+C=(V_1^0)^T\hat
N_0^1$.
\end{itemize}

\subsection{Two algebraic consequences of the tangency configuration}

Note that $x_z(u,v)$ for (\ref{eq:Qhyp}) is an affine image of the
equilateral hyperboloid with one sheet
$H(u,v):=\frac{(1-v^2)e_1+(1+v^2)e_2+2ve_3}{u-v}+\pa_v\frac{(1-v^2)e_1+(1+v^2)e_2+2ve_3}{2}=\\
\frac{(1-v^2)e_1+(1+v^2)e_2+2ve_3}{u-v}+H(\infty,v)=-H(v,u)$.

Let $\mathcal{B}:=(u-v)^2$ for (\ref{eq:Qhyp}) and
$\mathcal{B}:=1$ for (\ref{eq:Qpar}), $x_0^0:=x_0(u_0,v_0),\
x_0^1:=x_0(u_1,v_1),\ m_0^1:=\mathcal{B}_1\pa_{u_1}x_z^1\times
V_0^1$ a normal field of the distribution $\mathcal{D}^1$ of
facets $\mathcal{F}^1$ passing through $x_z^1$ and spanned by
$V_0^1,\ \pa_{u_1}x_z^1$ (and similarly
$m'^1_0:=\mathcal{B}_1\pa_{v_1}x_z^1\times V_0^1$ by considering
the other ruling family on $x_z^1$).

For $u_0,v_0,u_1,v_1$ independent variables
$\mathcal{B}_1\pa_{u_1}x_z^1$ depends only on $v_1$ (quadratically
for (\ref{eq:Qhyp}) and linearly for (\ref{eq:Qpar})), so
$\pa_{u_1}m_0^1=\mathcal{B}_1\pa_{u_1}x_z^1\times
\pa_{u_1}V_0^1=0$ and $m_0^1$ does not depend on $u_1$.

For (\ref{eq:Qhyp}) $\pa_{v_1}(\mathcal{B}_1\pa_{u_1}x_z^1)\times
x_z^1(\infty,v_1)=0$, so
$m_0^1=(\mathcal{B}_1\pa_{u_1}x_z^1)\times(x_z^1(\infty,v_1)-x_0^0)$
depends quadratically on $v_1$ (the coefficient of the highest
order term $v_1^3$ is $0$ and that of $v_1^2$ contains
$-\frac{1}{2}\pa_{v_1}^2(\mathcal{B}_1\pa_{u_1}x_z^1)\times
x_0^0$).

For (\ref{eq:Qpar})
$m_0^1=(\mathcal{B}_1\pa_{u_1}x_z^1)\times(x_z^1(0,v_1)-x_0^0)$;
since
$\pa_{v_1}(\mathcal{B}_1\pa_{u_1}x_z^1)\times\pa_{v_1}x_z^1(0,v_1)\neq
0$, we conclude that $m_0^1$ depends quadratically on $v_1$.

Thus we conclude that in both cases $m_0^1$ depends only on
$u_0,v_0,v_1$ and quadratically in $v_1$; this will make the
integrability condition (the differential equation subjacent to
the B\"{a}cklund transformation) a Ricatti equation in $v_1$.

Henceforth consider only the tangency configuration
$(V_0^1)^TN_0^0=0$ (from the proof of Bianchi I this will impose a
functional relationship among $u_0,v_0,u_1,v_1$ separately linear
in each variable, so a homography is established between them).

If we choose the rulings $w_0^0:=\pa_{u_0}x_0^0,\
w_0^1:=\pa_{u_1}x_0^1$ respectively at $x_0^0,\ x_0^1$, then we
get a rigid motion $(R_0^1,t_0^1)$ provided by the Ivory affinity.
If we change the ruling family on $x_0^1$, then the action of the
new rigid motion on the facet $T_{x_0^0}x_0$ does not change, so
its new rotation must be the old rotation composed with a
reflection in $T_{x_0^0}x_0$, because of which the facets
$\mathcal{F}^1,\ \mathcal{F}'^1$ reflect in $T_{x_0^0}x_0$ (thus
the distributions $\mathcal{D}^1,\ \mathcal{D}'^1$ reflect in
$Tx_0$):
\begin{eqnarray}\label{eq:Qref}
(\pa_{v_1}x_z^1)^T(I_3-2N_0^0(N_0^0)^T)m_0^1=0,
\end{eqnarray}
and $\pa_{v_1}x_0^1=R_0^1(I_3-2N_0^0(N_0^0)^T)\pa_{v_1}x_z^1$;
multiplying this on the left by $(\pa_{u_1}x_0^1)^T$ and using the
preservation of lengths of rulings under the Ivory affinity we get
\begin{eqnarray}\label{eq:Qdir}
4(\pa_{u_1}x_z^1)^TN_0^0(N_0^0)^T\pa_{v_1}x_z^1du_1dv_1=
|dx_z^1|^2-|dx_0^1|^2=-\frac{4z}{\mathcal{B}_1}du_1dv_1.
\end{eqnarray}
Thus we have the next result, essentially due to Bianchi (he uses
equivalent computations; in fact most relevant consequences of the
tangency configuration are either equivalent to it or to it
composed with simple symmetries):

\begin{lemma}
If $x_z^1\in T_{x_0^0}x_0$, then:

I (Factorization) The change in the linear element from $x_z^1$ to
$x_0^1$ is four times the product of the orthogonal projections of
the differentials of the rulings of $x_z^1$ on the normal of $x_0$
at $x_0^0$.

II (Reflection) The facets at $x_z^1$ spanned by $V_0^1$ and one
of the rulings of $x_z^1$ reflect in $T_{x_0^0}x_0$; therefore the
distributions $\mathcal{D}^1,\ \mathcal{D'}^1$ reflect in
$Tx_0^0$.
\end{lemma}

\begin{remark}
Although the equal inclination of facets to tangent planes of
seeds from {\bf Remark \ref{rm:refl}} is not preserved when one
considers general quadrics instead of pseudo-spheres, the
reflection property of facets in tangent planes of seeds remains
valid and it is explained by the existence of the rigid motion
provided by the Ivory affinity regardless of the choice of
rulings. Note however that although this explanation is good
enough from an analytic point of view, the choice of ruling still
lacks geometric motivation and {\bf Remark \ref{rm:refls}}
provides it.
\end{remark}

The algebraic relation
\begin{eqnarray}\label{eq:Qintalg}
(N_0^0)^T(2zm_0^1+m_0^1\times\pa_{v_1}m_0^1)=0
\end{eqnarray}
will appear as the total integrability condition. Using
(\ref{eq:Qref}), (\ref{eq:Qdir}) this becomes:
$0=\frac{z(m_0^1)^T\pa_{v_1}x_z^1}{(N_0^0)^T\pa_{v_1}x_z^1}-
\mathcal{B}_1(\pa_{u_1}x_z^1)^TN_0^0(V_0^1)^T
\pa_{v_1}m_0^1=\frac{z(V_0^1)^T(\mathcal{B}_1\pa_{v_1}x_z^1\times
\pa_{u_1}x_z^1+\pa_{v_1}(\mathcal{B}_1\pa_{u_1}x_z^1\times
V_0^1))}{(N_0^0)^T\pa_{v_1}x_z^1}$, which is straightforward.
Replacing $(m_0^1,v_1)$ with $(m'^1_0,u_1)$ we get a similar
relation.

\subsection{Rolling quadrics and distributions}

Let $(R_0,t_0)(x_0^0,dx_0^0)=(x^0,dx^0)$ be the rolling of the
piece of quadric $x_0^0=x_0(u_0,v_0)$ on the isometric surface
(seed) $x^0\subset\mathbb{R}^3$. The facets of the rolled
distribution $(R_0,t_0)\mathcal{D}^1$ will become tangent planes
to leaves $x^1:=(R_0,t_0)x_z^1=(R_0,x^0)V_0^1$ iff the
integrability condition $0=(R_0m_0^1)^Tdx^1$ holds. We have
$R_0^{-1}dx_1=d(x_0^0+V_0^1)+R_0^{-1}dR_0V_0^1=dx_z^1+\omega_0\times
V_0^1$. But $(\omega_0)^{\bot}=0$ and
$dx_z^1=\pa_{u_1}x_z^1du_1+\pa_{v_1}x_z^1dv_1$, so the
integrability condition becomes $-(V_0^1)^T\omega_0\times
N_0^0(m_0^1)^TN_0^0+(m_0^1)^T\pa_{v_1}x_z^1dv_1=0$; using
(\ref{eq:Qref}) this becomes $-(V_0^1)^T\omega_0\times
N_0^0+2(N_0^0)^T\pa_{v_1}x_z^1dv_1=0$; multiplying it by
$\mathcal{B}_1(N_0^0)^T\pa_{u_1}x_z^1$, using (\ref{eq:Qdir}) and
$-\mathcal{B}_1(N_0^0)^T\pa_{u_1}x_z^1V_0^1=\mathcal{B}_1(V_0^1\times
\pa_{u_1}x_z^1)\times N_0^0=-m_0^1\times N_0^1$ we finally get the
Ricatti equation:
\begin{eqnarray}\label{eq:Qric}
(m_0^1)^T\omega_0+2zdv_1=0.
\end{eqnarray}
We have $dm_0^1=\pa_{v_1}m_0^1dv_1+\mathcal{B}_1dx_0^0\times
\pa_{u_1}x_z^1$, so
$(dm_0^1)^T\wedge\omega_0=dv_1\wedge(\pa_{v_1}m_0^1)^T\omega_0$;
imposing the total integrability condition $d$ on (\ref{eq:Qric})
and using the equation itself we need
$-(m_0^1)^T\omega_0\wedge(\pa_{v_1}m_0^1)^T\omega_0+2z(m_0^1)^Td\omega_0=0$,
or, using (\ref{eq:oma}) and (\ref{eq:Qom}):
$(N_0^0)^T(2zm_0^1+m_0^1\times
\pa_{v_1}m_0^1)(N_0^0)^T\omega_0\times\wedge\omega_0=0$; thus the
total integrability is equivalent to (\ref{eq:Qintalg}).

\section{ The link to {\it The Method} of Archimedes}

If we roll $x_0^0$ on different sides of the seed $x^0$, then we
get the B\"{a}cklund transformation for the other ruling family,
so the rolled distributions reflect in the bundle of tangent
planes of the seed $x^0$. Thus  (\ref{eq:Qref}) is obtained if one
makes the ansatz $x_0^0=x^0$; the same ansatz for the isometric
correspondence provided by the Ivory affinity and the inversion of
the B\"{a}cklund transformation (two focal surfaces of a line
congruence are in a symmetric relationship) implies Bianchi's
result about the existence of $(R_0^1,t_0^1)$ and the symmetry of
the tangency configuration; now (\ref{eq:Qdir}) is obtained as
previously described.

\begin{remark}
While {\it 'first'} in {\it The Method} clearly can be linked to
the Bianchi-Lie ansatz, just by fortuitous chance or by
Archimedes' clairvoyance (or a combination thereof) the {\it
'mechanical method'} meant by Archimedes for the use of balance
and slicing corresponds to rolling (clearly a mechanical method)
with its inherent indeterminacy and {\it 'although they had to be
proved by geometry afterwards because their investigation by the
said method did not furnish an actual proof'} used by Archimedes
for the double reduction ad absurdum corresponds to the geometric
arguments involved in the discussion on the rigid motion provided
by the Ivory affinity (which may also correspond to Archimedes'
geometric identities at the infinitesimal level and using the
moments of the balance).
\end{remark}

Note that the method at the level of points of facets was known to
Bianchi and Lie; however, they had never used the full method, at
the level of the planes of the facets too (for this reason the
{\it 'of course easier'} ingredient is missing from Bianchi's
proofs). Thus if one assumes Theorem I of Bianchi's theory of
deformations of quadrics a-priori to be true and to be the
metric-projective generalization of Lie's approach, then one
naturally {\it geometrically} gets the necessary algebraic
identities needed to prove Theorem I.

For surfaces this method is just a fancy way of reformulating
already known identities and which appear naturally enough at the
analytic level. But keep in mind that the tangency configuration,
(\ref{eq:Qref}) and (\ref{eq:Qdir}) are equivalent from an
analytic point of view and this is not the case in higher
dimensions: thus it is very difficult to find the necessary
algebraic identities of the static picture from an analytic point
of view. Therefore {\it The Method} of Archimedes, due to its
geometric naturalness and the fact that it contains more
information, is useful in the study of higher dimensional
problems.

Note that although Archimedes and Bianchi-Lie have dealt with
different problems, their approach was the same: $P\Rightarrow Q$
with $P$ being either {\it 'The area of a segment of a parabola is
an infinite sum of areas of lines'} or {\it 'A line is a
deformation of a quadric'}. Such sentences $P$ were not fully
accepted as true according to the standard of proof of the times,
but they had valid relevant consequences $Q$ which elegantly
solved problems not solvable with other methods of those times.

Note that {\it The Method} of Archimedes went a little closer to
the B\"{a}cklund transformation: in an a-priori intuitive
elementary non-rigorous geometric argument he transferred all
lines (slices of the segment of the parabola) with their centers
at the left end of the balance (thus with $\infty$ multiplicity,
similarly to the original position of facets in the Bianchi-Lie
ansatz). Thus the quadratic mass of a slice of a parabola segment
at the left hand side of the balance factorizes in the product of
the linear mass of a slice of a triangle and the linear length of
the leg of the right hand side of the balance; by integration the
mass of the parabola placed with its center at the left hand side
of the balance remains in equilibrium with the mass of the
triangle placed with its center at the right hand side of the
balance.

Note that the facets of the $3$-dimensional rolled distribution
are differently re-distributed into $2$-dimensional families of
facets as tangent planes to leaves when the shape of the seed
changes (for Bianchi's complementary transformation they remain
the same), but principles and properties independent of the shape
of the seed remain valid for facets even in the singular picture:
this is Archimedes' contribution to the Bianchi-Lie ansatz and
allows us to call this full Bianchi-Lie ansatz {\it the
Archimedes-Bianchi-Lie method}.

Thus one can conclude that {\it Archimedes' balance} is one and
the same with {\it Bianchi's B\"{a}cklund transformation for
quadrics} as principles of a general nature: they are valid at the
infinitesimal level and survive integration and conversely, being
principles of a general nature both induce by differentiation and
by particular singular configurations the infinitesimal picture
where the simplest explanation of these principles reveals itself.

\begin{remark}
We have mostly used 'tangent planes' instead of 'tangent spaces',
since a point, a line or a surface in $\mathbb{R}^3$ come with a
$2$-dimensional family of tangent planes and either can appear as
leaves of integrable distributions. Most of the general statements
remain valid when the dimension of the leaves collapses (some of
'surfaces' may have to be replaced with 'lines').
\end{remark}

\begin{remark}
Note that in \cite{A0} Archimedes quotes an even earlier result of
Democritos related to the volume of the cone as an analogy and
inspiration to his method; thus it may be the case that Archimedes
followed the same footsteps: from a finite law obtained by
empirical observation one gets an infinitesimal law by applying
the same finite law to thinner finite objects and a collapsing end
process; finally with all relevant information recorded by the
infinitesimal objects and which is easier to prove one rigorously
proves the general conjectured finite law. This is an inverting
point of view similar to Lie's.
\end{remark}

\section*{Acknowledgements}

The research partially appearing in this paper was done during a
graduate program at the University of Notre Dame du Lac; the
author wishes to thank for the academic support during this
graduate program (which included also Summers) and useful
discussions and advice from Advisor Professor Brian Smyth. Also
the author wishes to thank for academic support from the
Mathematics Department of Bucharest University in a graduate
program beginning with September 2007 and to my former (late)
Advisor Professor Stere Ianu\c{s} for advice and feedback on
drawing an early form of this paper.

\section*{Note added}

According to one of our notes {\it On Bianchi's B\"{a}cklund
transformation for quadrics}\footnote{See \cite{B3} and
{\href{http://arxiv.org/abs/1110.5474} {On Bianchi's B\"{a}cklund
transformation for quadrics}}} we have improved these results and
some of them still fall into the Archimedes-Bianchi-Lie's approach
of collapsing ansatz of leaves of the $3$-dimensional integrable
rolling distribution of facets (the defining surface is quadric
and the leaves are rulings on a confocal quadric, which includes
in a general sense (isotropic) plane of a pencil of (isotropic)
planes for quadrics of revolution or Darboux quadrics).

If facets of a $3$-dimensional integrable rolling distribution of
facets are centered on tangent planes of a surface $x_0$ and
further pass through the origin of these tangent planes (thus to
each tangent plane we associate an $1$-dimensional family of
facets), then the Weingarten congruence property is satisfied
without further requirements.

For a $3$-dimensional integrable tangential rolling distribution
of facets with collapsing ansatz of leaves to curves (facets are
centered on tangent planes of the defining surface $x_0$ and on an
auxiliary surface $x_z$; thus the original leaves are two families
of curves which induce special coordinates on $x_z$ and the
reflection property plays again an important r\^{o}le) we get the
symmetric tangency configuration (facets further pass through the
origin of these tangent planes) and isometric correspondence of
leaves of a general surface (independent of the shape of the
surface $x^0$ isometric to $x_0^0\subset x_0$); thus we get the
most general form of a a theory of isometric deformations of
surfaces via B\"{a}cklund transformation with defining surface.

However, we have found no other defining surfaces besides quadrics
and we have proved that there are no other surfaces, thus
Bianchi's generalization of the B\"{a}cklund transformation of the
pseudo-sphere to quadrics is the completion of the theory of
isometric deformations of surfaces via B\"{a}cklund transformation
with defining surface.

\bibliographystyle{amsplain}

\end{document}